\documentclass[11pt,a4paper,reqno]{amsart}
\usepackage{lipsum}

\usepackage[utf8]{inputenc}
\usepackage{amsmath}
\usepackage{amsthm}
\usepackage{amsfonts}
\usepackage{pgf,tikz}
\usepackage{tikz-cd}
\usepackage{amssymb}
\usepackage{dirtytalk}
\usepackage{float}
\usepackage{mathtools}
\usepackage{imakeidx}
\tikzset{every picture/.style={line width=0.75pt}}

\usepackage{graphicx}
\usepackage{xcolor}
\usepackage{etoolbox}

\definecolor{rojo}{rgb}{208, 2, 27}
\definecolor{verde}{RGB}{65,117,5}
\definecolor{azul}{RGB}{74,144,226}

\newtheorem{theorem}{Theorem}
\numberwithin{theorem}{section} 
\newtheorem{lemma}[theorem]{Lemma}

\newtheorem{proposition}[theorem]{Proposition}

\newtheorem{definition}[theorem]{Definition}

\theoremstyle{remark}
\newtheorem{remark}[theorem]{Remark}
\newtheorem{convention}[theorem]{Convention}
\newtheorem{example}[theorem]{Example}
\usepackage{color}

\newcommand{\pc}[1] {\textcolor{purple}{#1}}

\usepackage[colorlinks=true,linkcolor=blue,citecolor=blue
]{hyperref}%
\usepackage{booktabs}

\usepackage{caption} 
\newcommand{\discr}{\operatorname{discr}}

\DeclareMathOperator{\rk}{\operatorname{rk}}
\DeclareMathOperator{\iso}{\cong}

\DeclareMathOperator{\NS}{\operatorname{NS}}

\newcommand{\Ee}{\mathcal{E}}
\usepackage{adjustbox}

\usepackage{enumerate}

\usepackage{todonotes}
\usepackage[normalem]{ulem}

\title{On strictly elliptic K3 surfaces and del Pezzo surfaces}

\author[P. Comparin]{Paola Comparin}
\address{Departamento de Matemática y Estadística, Universidad de la Frontera, Temuco, Chile}  \email{paola.comparin@ufrontera.cl}

\author[P. Montero]{Pedro Montero}
\address{Departamento de Matem\'aticas, Universidad T\'ecnica
  Fe\-de\-ri\-co San\-ta Ma\-r\'\i a, Valpara\'\i
  so, Chile}  \email{pedro.montero@usm.cl}

\author[Y. Prieto--Monta\~{n}ez]{Yulieth Prieto--Monta\~{n}ez} %
\address{The Abdus Salam International Centre for Theoretical Physics, Str. Costiera, 11, 34151 Trieste TS, Italy} %
\email{yprieto@ictp.it}

\author[S. Troncoso]{Sergio Troncoso}
\address{Departamento de Matem\'aticas, Universidad T\'ecnica
  Fe\-de\-ri\-co San\-ta Ma\-r\'\i a, Valpara\'\i
  so, Chile}  
\email{sergio.troncosoi@usm.cl, stroncosoi11@gmail.com}

\date{}
\makeindex

\begin{document}

\maketitle
\begin{abstract}

This article primarily aims at classifying, on certain K3 surfaces, the elliptic fibrations induced by conic bundles on smooth del Pezzo surfaces. The key geometric tool employed is the Alexeev-Nikulin correspondence between del Pezzo surfaces with log-terminal singularities of Gorenstein index two and K3 surfaces with non-symplectic involutions of elliptic type: the latter surfaces are realized as appropriate double covers obtained from the former ones. The main application of this correspondence is in the study of linear systems that induce elliptic fibrations on K3 surfaces admitting a strictly elliptic non-symplectic involution, i.e., whose fixed locus consists of a single curve of genus $g\geq 2$. The obtained results are similar to those achieved by Garbagnati and Salgado for jacobian elliptic fibrations.
\end{abstract}

\tableofcontents

We work over the field of complex numbers $\mathbf{C}$.

\section{Introduction}

One of the principal algebraic invariants of a projective algebraic variety $X$ is its group of biregular automorphisms, denoted as $\operatorname{Aut}(X)$. In many cases, this group can be used to discover geometric properties of the underlying algebraic variety and its projective models. Furthermore, as we will explore, the mere existence of suitable involutions has non-trivial consequences regarding special linear systems on the corresponding variety.

One of the main successful approaches to study the automorphism group of a smooth projective variety $X$ is to consider the natural induced action on cohomology lattices. Indeed, a classical result by Lieberman \cite{Lie78} states that the neutral component $\operatorname{Aut}^\circ(X)$ has finite index in the kernel of the linear representation  $\operatorname{Aut}(X)\to \operatorname{GL}(H^2(X,\mathbf{Z}))$, i.e., $\operatorname{Aut}(X)$ splits into its neutral component $\operatorname{Aut}^\circ(X)$ and its discrete image in $\operatorname{GL}(H^2(X,\mathbf{Z}))$ (see e.g. \cite{Can18} for further details).

The aforementioned strategy has been classically employed to investigate the automorphism groups of smooth \emph{del Pezzo surfaces}, which are projective algebraic surfaces $Z$ such that the anti-canonical divisor $-K_Z$ is ample. In cases where these groups are finite, they are realized as subgroups of the Weyl group of an appropriate lattice (see \cite[Chapters 8--9]{Dol12} for details and a historical account). Another remarkable class of algebraic varieties for which this method is commonly used are \emph{projective K3 surfaces}, i.e., simply connected smooth projective surfaces $X$ with trivial canonical bundle $K_X\sim 0$. Indeed, for these surfaces the relevance of the study of automorphisms is a consequence of the celebrated Torelli theorem due to Pyatetskii-Shapiro and Shafarevich \cite{SS71}, and general lattice theoretic results by Nikulin \cite{Nik76,Nik79,Nik81}.

Remarkably, these two classes of surfaces are naturally related thanks to the work of Alexeev and Nikulin \cite{AleNik06} (see also \cite{Nak07,Zha98})  on the classification of del Pezzo surfaces $Z$ of Gorenstein index 2, i.e., the Weil divisor $-2K_Z$ is an ample Cartier divisor. More precisely, let us recall that in any K3 surface $X$ there is a holomorphic and non-vanishing (i.e., symplectic) 2-form $\omega_X$ such that $\operatorname{H}^0(X,\Omega_X^2)=\mathbf{C}\cdot \omega_X$, and thus we say that $\sigma \in \operatorname{Aut}(X)$ is \emph{symplectic} if $\sigma^\ast(\omega_X)=\omega_X$, and that is \emph{non-symplectic} otherwise. In this terms, the correspondence \cite{AleNik06} between K3 surfaces $X$ and (possibly) singular del Pezzo surfaces $Z$ of Gorenstein index 2 goes as follows: 
the so-called Smooth Divisor Theorem (see \cite[Theorem 1.5]{AleNik06}) ensures that, given such a del Pezzo surface, there is a smooth irreducible curve $C$ in the linear system $|-2K_Z|$. Subsequently, a double cover $W\to Z$ can be constructed, branched over both $C$ and $\operatorname{Sing}(Z)$. 
The crucial observation lies in the fact that the minimal resolution $X\to W$ is a K3 surface, and the involution $\iota:X\to X$ associated with the covering is such that $\iota^\ast \omega_X = - \omega_X$, i.e., $\iota$ is a \emph{non-symplectic involution}. 
Through an analysis of $X$ and the quotient surface $X\slash \langle \iota \rangle$, the authors in \cite{AleNik06} reduce the classification problem for $Z$ to the study of K3 surfaces with non-symplectic involutions.

It is worth mentioning that the above construction is a vast generalization of the fact that the double cover of $\mathbf{P}^2$ branched over a \emph{smooth} sextic curve $C$ is a K3 surface with a non-symplectic involution. This classical construction (revisited by Dolgachev \cite{Dol73} and Reid \cite{Reid} using modern methods) was considered by Enriques and Campedelli (see \cite{Cam40}), who investigated double coverings of $\mathbf{P}^2$ that are birational to a K3 surface. Another remarkable recent construction in birational geometry, closely related to our context, is the work of Peters and Sterk \cite{PS20} where they consider nodal Enriques surfaces constructed from a K3 surface obtained as a double cover of a smooth del Pezzo surface of degree 6.

Even in the simplest case where the quotient del Pezzo surfaces are smooth (which will be our case of interest), the correspondence by Alexeev and Nikulin enables a connection between the geometric properties of the del Pezzo surface $Z=X\slash \langle \iota \rangle$ and the lattice-theoretic properties of the invariant lattice $\operatorname{NS}(X)^{\iota^\ast}$ associated with the non-symplectic involution $\iota$ 
on the K3 surface (see \S 2.2 for details). For instance, it follows from results by Nikulin in \cite{Nik83} that the topology of the fixed locus $X^\iota$ is determined by suitable discrete invariants $(r,a,\delta)$ of the lattice $\operatorname{NS}(X)^{\iota^\ast}$ which, in turn, can be interpreted through $Z$ (e.g., the genus of the fixed curve $X^\iota=C_g$ will be given by $g=K_Z^2+1$).

After understanding the topology of the fixed locus, a natural next step is to describe the possible linear systems on these K3 surfaces, following the work initiated by Saint-Donat in \cite{SD74}. Notably, linear systems inducing elliptic fibrations (see Definition \ref{def:elliptic fibration}, and note that we do not require the existence of a section) are of special interest, as they can be characterized numerically thanks to the results in \cite[\S 3]{SS71} and since they have important arithmetic applications (see e.g. \cite{Sch09}). 
Significant progress has been made in this direction, particularly in the case where the K3 surface is \emph{generic} among those admitting a non-symplectic automorphism with a given fixed locus (see Convention \ref{convention:generic}). More precisely, in \cite{Nik83, Ogu89, Klo06, CG, GarSal19, GarSal20} the authors classify elliptic fibrations on K3 surfaces with a non-symplectic involution in many cases, using the fact that the fixed locus of the involution is either empty, the disjoint union of two elliptic curves or contains at least a rational curve (we refer the reader to Remark \ref{remark:remaining-cases} for further details). 
In order to complete the classification of elliptic fibrations on such generic surfaces, it remains to study the case when $X^\iota = C_g$ consists of a single smooth irreducible curve of genus $g\geq 2$. In regard of the previous discussion, 

\begin{quote}
    the main purpose of this article is to address, through the Alexeev-Nikulin correspondence, the remaining case in the classification of (not necessarily jacobian) elliptic fibrations on generic K3 surfaces (in the sense of Convention \ref{convention:generic}) that admit a non-symplectic involution $\iota:X\to X$.
\end{quote}

To achieve this, we will restrict ourselves to the case where $X^\iota = C_g$ is a smooth irreducible curve of genus $g\geq 2$, and we will say that $\iota:X\to X$ is a \emph{strictly elliptic involution} (see Definition \ref{def: strictly elliptic type}). As we will observe in Proposition \ref{prop: K3 of strictly elliptic type}, the main feature of the pair $(X,\iota)$ is that the Alexeev-Nikulin correspondence results in a \emph{smooth} del Pezzo surface $Z=X\slash \langle \iota \rangle$. In this context, our main result (see Theorem \ref{thm_main}) establishes that the quotient projection $\pi:X\to Z=X\slash \langle \iota \rangle$ induces a correspondence between elliptic fibrations $\mathcal{E}:X\to \mathbf{P}^1$ and conic bundles (see Definition \ref{Def: Conic Bundle}) $f:Z\to \mathbf{P}^1$, and moreover $\mathcal{E}=f\circ \pi$. It is worth noting that analogous results have been obtained for other K3 surfaces with non-symplectic involutions, provided the fixed locus $X^\iota$ contains rational curves (see e.g. \cite[\S 5]{GarSal19}, \cite[\S 7.1]{GarSal20} and \cite[\S 5.10]{CG}). In contrast, the main advantage of our approach using del Pezzo surfaces, as opposed to previous works relying on lattice-theoretic methods (where the existence of a section of the elliptic fibration or the presence of a rational curve in the fixed locus $X^\iota$ is important), is its compatibility with standard tools from Mori theory (see e.g. \cite[Chapter 6]{Deb01}). Notably, in our case, we can classify the effective classes of curves inducing conic bundles on the surface $Z$ (see Proposition \ref{Classes of conic bundles on Z}), following a similar approach as in the case of $(-1)$-curves in \cite[\S 26]{Man86}, and use them to describe all admissible singular fibers of the induced elliptic fibrations on the corresponding K3 surface $X$ (see \S 5).

Finally, it is noteworthy that, as a consequence of recent work \cite{CliMal23} by Clingher and Malmendier, the considered elliptic fibrations are not jacobian, i.e., they do not admit sections. However, it is not difficult to observe (see \S 5) that they admit bisections which can be induced by the $(-1)$-curves in the associated del Pezzo surface (see Example \ref{example:bisections}). 
 Despite the absence of jacobian elliptic fibrations, and consequently the inability to consider Weierstrass models, 
 these K3 surfaces are quite special due to the fact that by
\cite[\S 2.8]{AleNik06}  they have \emph{finite automorphism groups}. Remarkably, they fall within the recent work by Roulleau \cite{Rou}, where some explicit projective models are studied and where the full lattice $\operatorname{NS}(X)$ is described (see \S 4).

\subsection*{Acknowledgements} The authors would like to sincerely thank \textsc{Michela Artebani}, \textsc{Giacomo Mezzedemi}, \textsc{Cec\'ilia Salgado}, and \textsc{Giancarlo Urz\'ua} for various fruitful discussions regarding the constructions used in this article. Special thanks go to \textsc{Cinzia Casagrande}, who informed us about reference \cite{Der06} and allowed us to correct errors in the counting of conic classes at low degree. P.C. has been partially funded by Universidad de La Frontera, Proyecto DIM23-0001 and Fondecyt ANID Project 1240360. P.M. has been partially funded by Fondecyt ANID Projects  1231214 and 1240101. S.T. was partially supported by Fondecyt ANID Project 3210518. The authors express their gratitude to \textsc{Cec\'ilia Salgado} for her generous hospitality extended to S.T. during his visit to the Bernoulli Institute at the University of Groningen.

\section{Background and preliminaries}
\subsection{Conic bundles and smooth del Pezzo surfaces}

Let us recall that a smooth projective surface $Z$ is called a del Pezzo surface if the anti-canonical divisor $-K_Z$ is ample. The positive integer $d(Z)=(-K_Z)^2$ is called the degree of $Z$, and it is the main invariant that allows for their classification. More precisely, we have the following classical result (see e.g. \cite[\S 24]{Man86} and \cite[Proposition 8.1.25]{Dol12}).

\begin{theorem}\label{thm:del_Pezzo}
    Let $Z$ be a smooth del Pezzo surface of degree $d$. Then, $1\leq d \leq 9$ and we have that:
    \begin{itemize}
        \item[(i)] If $d=9$, then $Z\simeq \mathbf{P}^2$.
        \item[(ii)] If $d=8$, then $Z$ is isomorphic to either $\mathbf{P}^1\times \mathbf{P}^1$ or to the blow-up of $\mathbf{P}^2$ at one point \emph{(}i.e., the Hirzebruch surface $\mathbf{F}_1$\emph{)}.
        \item[(iii)] If $1\leq d\leq 7$, then $Z\simeq \operatorname{Bl}_{p_1,\ldots,p_r}(\mathbf{P}^2)$ where $r=9-d$ and where the points $p_1,\ldots,p_r\in \mathbf{P}^2$ are in \emph{general position}.
    \end{itemize}
    Here, we say that the points are in \emph{general position} if the following hold:
    \begin{enumerate}
        \item no three points are on a line;
        \item no six points are on a conic;
        \item no nodal or cuspidal cubic passes through eight points with one of them being the singular point.
    \end{enumerate}
    Conversely, any blow-up of $r\leq 8$ points in general position is a smooth del Pezzo surface.
\end{theorem}

\begin{convention}\label{convention:Z_d}
    We will denote by $Z$ an arbitrary smooth del Pezzo surface. Additionally, we will denote by $Z_d$ the del Pezzo surface of degree $d\in \{1,\ldots,8\}$ obtained as the blow-up of $9-d$ points in general position in $\mathbf{P}^2$, and by $E_1,\ldots,E_{9-d}$ the corresponding
    exceptional divisors.
\end{convention}

It is worth mentioning that the geometry of exceptional curves on del Pezzo surfaces is completely understood. For instance, it is known that every irreducible curve with negative self-intersection on a smooth del Pezzo surface $Z$ is a $(-1)$-curve (see e.g. \cite[Theorem 24.3]{Man86}). More precisely, we have the following result (see \cite[\S 26]{Man86} and \cite[\S 8.2.6]{Dol12}).

\begin{theorem}\label{thm:(-1)-curves}
    Let $Z_d$ be a smooth del Pezzo surface of degree $1\leq d \leq 8$, and let $\varepsilon:Z_d \to \mathbf{P}^2$ be its representation as the blow-up of $r=9-d$ points $p_1,\ldots,p_r\in \mathbf{P}^2$ in general position. Let $\Gamma\subseteq Z_d$ be a $(-1)$-curve, then the image $\varepsilon(\Gamma)\subseteq \mathbf{P}^2$ is of one of the following types:
    \begin{enumerate}
        \item one of the points $p_i$;
        \item a line passing through 2 of the points $p_i$;
        \item a conic passing through 5 of the points $p_i$;
        \item a cubic passing through 7 of the points $p_i$ such that 1 of them is a double point;
        \item a quartic passing through 8 of the points $p_i$ such that 3 of them are double points;
        \item a quintic passing through 8 of the points $p_i$ such that $6$ of them are double points;
        \item a sextic passing through 8 of the points $p_i$ such that 7 of them are double points and one is a triple point.
    \end{enumerate}
    Moreover, the number $n$ of $(-1)$-curves on $Z_d$ is given by the following table
    \[
\begin{tabular}{|c||cccccccc|} \hline
$d$ & $8$ & $7$ & $6$ & $5$ & $4$ & $3$ & $2$ & $1$ \\ \hline
$r$ & $1$ & $2$ & $3$ & $4$ & $5$ & $6$ & $7$ & $8$ \\ \hline
$n$ & $1$ & $3$ & $6$ & $10$ & $16$ & $27$ & $56$ & $240$ \\ \hline
\end{tabular}
    \]
\end{theorem}

Following the same line of ideas that allow classifying the images of $(-1)$-curves in the above result, we can describe the possible conic bundles on smooth del Pezzo surfaces. For the reader's benefit, we recall the relevant notions about conic bundles below.

\begin{definition}{\cite[\S 1]{Sarkisov}}\label{Def: Conic Bundle}
A conic bundle on a smooth projective surface $Z$ is a surjective morphism onto a smooth curve $f:Z\to C$ whose general fiber is a smooth, irreducible curve of genus $0$.
\end{definition}

\begin{remark}
    According to L\"uroth's Theorem, if $Z$ is a rational surface (e.g. a smooth del Pezzo surface) then the curve $C$ must be isomorphic to $\mathbf{P}^1$.
\end{remark}

\begin{definition}\label{Canonical Class} Let $Z$ be a smooth projective surface, we say that
an element $[D]\in \operatorname{NS}(Z)$ is a conic class  if $D$ is nef, $D^2=0$, and $D\cdot K_Z=-2$.
    
\end{definition}

\begin{lemma}\label{Numerical characterization of conic bundle}Let $Z$ be smooth rational surface, then
    there exists a correspondence between the set of conic bundles on $Z$ and the set of conic classes of $\operatorname{NS}(Z)$.
\end{lemma}
\begin{proof}
First, let $[D]\in \operatorname{NS}(Z)$ be a conic class, then by Riemman-Roch and \cite[Theo. III.1.(a)]{Anticanonical} the associated map $\phi_{|D|}:Z\to\mathbf{P}^1$ is a well-defined morphism, and since every fiber is linearly equivalent to $D$ it is a conic bundle on $Z$.
On the other hand, let $f:Z\to \mathbf{P}^1$ be a conic bundle, then it is clear that the class of any fiber of $f$, say $F$, is a conic class.
\end{proof}

\begin{remark}\label{delta}
The number of singular fibers of a conic bundle is a numerical invariant. For instance, for a conic bundle $f:Z_d\to \mathbf{P}^1$ on the del Pezzo surface $Z_d$, it is well-known that the number of singular fibers is 
$8-K_{Z_d}^2=8-d$, as seen in \cite[\S 1]{KollarMella}.   
\end{remark}

 \begin{proposition}\label{Classes of conic bundles on Z}
     Let $Z_d$ be a del Pezzo surface of degree $d\leq 8$ obtained as the blow-up of $9-d$ points in general position in $\mathbf{P}^2$. Then the conic classes are listed in Table \ref{tab_canonical_bundles}.  Moreover, the number $N$ of conic bundles on $Z_d$ is given by the following table
 \[
\begin{tabular}{|c||cccccccc|} \hline
$d$ & $8$& $7$ & $6$ & $5$ & $4$ & $3$ & $2$ & $1$ \\ \hline
$r$ & $1$& $2$ & $3$ & $4$ & $5$ & $6$ & $7$ & $8$ \\ \hline
$N$ & $1$&$2$ & $3$ & $5$ & $10$& $27$ & $126$ &  $2160$ \\ \hline
\end{tabular}
    \] 
\begin{table}[h]
\resizebox{\textwidth}{!}{%
\begin{tabular}{|l|l|l|}
\hline
$d=8$ & $D=L-E_1$                                     \\ \hline
$d=7$ & $D=L-E_i$, $i=1,2$                            \\ \hline
$d=6$ & $D=L-E_i$, $i=1,2,3$                          \\ \hline
$d=5$ & $\displaystyle D=L-E_i$, $i=1,2,3,4$          \\ \hline
     & $\displaystyle D=2L-\sum_{k=1}^{4}E_k$         \\ \hline
$d=4$ & $\displaystyle D=L-E_i$, $i=1,2,3,4,5$                                                                                                                                      \\ \hline
      & $\displaystyle D=2L-\sum_{j\in J}E_{j}$,  $J\subset \{1,\ldots,5\}, |J|=4$                                                                                                \\ \hline
$d=3$ & $\displaystyle D=L-E_i$, $i=1,2,3,4,5,6$.                                                                                                                                      \\ \hline
      & $\displaystyle D=2L-\sum_{j\in J}E_{j}$,  $J\subset \{1,\ldots,6\}, |J|=4$                                                                                                 \\ \hline
      & $\displaystyle D=3L-2E_{i}-\sum_{j\in J}E_{j}$, $i\in \{ 1,\dots, 6\}$,  $J=\{ 1,\dots, 6\}\setminus \{i\}$                                                                            \\ \hline
$d=2$ & $\displaystyle D=L-E_i$, $i=1,\dots, 7$                                                                                                                                        \\ \hline
      & $\displaystyle D=2L-\sum_{j\in J}E_{j}$,  $J\subset \{1,\ldots, 7\}, |J|=4$                                                                                                \\ \hline
      & $\displaystyle D=3L-2E_{i}-\sum_{j\in J}E_{j}$, $i\in \{ 1,\dots, 7\}$,  $J=\{ 1,\dots, 7\}\setminus \{i\}, |J|=5$                                                                     \\ \hline
      & $\displaystyle D=4L-2\sum_{j\in J}E_{j}-\sum_{k\in K}E_{k}$,   $J\subset \{1,\ldots, 7\}, |J|=4$, $K=\{1,\ldots,7\}\setminus J$                                                 \\ \hline
      & $\displaystyle D=5L-E_{i}-2\sum_{j\in J}E_{j}$,   $i\in \{1,\dots,7\}$, $J=\{1,\ldots,7\}\setminus\{i\}$                                                                       \\ \hline
$d=1$ & $D=L-E_i$, $i=1,\dots,8$                                                                                                                                                       \\ \hline
      & $\displaystyle D=2L-\sum_{j\in J}E_{j}$,  $J\subset \{1,\ldots,8\}, |J|=4$                                                                                                 \\ \hline
      & $\displaystyle D=3L-2E_{i}-\sum_{j\in J}E_{j}$, $i\in \{ 1,\dots, 8\}$,  $J=\{ 1,\dots, 8\}\setminus \{i\}, |J|=5$                                                                     \\ \hline
      & $\displaystyle D=4L-2\sum_{j\in J}E_{j}-\sum_{k\in K}E_{k}$,   $J,K\subset \{1,\ldots,8\}, |J|=4, |K|=3, J\cap K=\emptyset$                                   \\ \hline

      &$\displaystyle D=4L-3E_i - \sum_{j\in J}E_j,\;i\in \{1,\ldots,8\},\;J=\{1,\ldots,8\}\setminus \{i\} $ \\\hline
      & $\displaystyle D=5L-3E_{i}-2\sum_{j\in J}E_{j}-\sum_{k\in K}E_{i_k}$, $i\in\{1,\ldots,8\}$, \\& $J, K\subset\{ 1,\dots, 8\}\setminus\{i\}, |J|=3, |K|=4, J\cap K=\emptyset$                                       \\ \hline
    &  $\displaystyle D=5L - E_i - 2 \sum_{j\in J} E_j,\; i\in \{1,\ldots,8\},\;J\subset \{1,\ldots,8\}\setminus \{i\},\;|J|=6 $\\ \hline
      
      & $\displaystyle D=6L-3\sum _{j\in J}E_{j}-2\sum_{k\in K}E_{k}-\sum_{r\in R}E_{r}$, \\& $J, K, R\subset\{1,\dots, 8\}, |J|=2, |K|=4, |R|=2, J\cap K= J\cap R=K\cap R=\emptyset$ \\ \hline

      & $\displaystyle D=7L - 4E_i - 3E_j - 2 \sum_{k\in K}E_k,\;i\in \{1,\ldots,8\},\;j\in \{1,\ldots,8\}\setminus \{i\},\;K = \{1,\ldots,8\}\setminus \{i,j\} $ \\  \hline
      & $\displaystyle D=7L - E_i - 2 \sum_{j\in J} E_j - 3 \sum_{k\in K} E_k,\;i\in \{1,\ldots,8\},\; J,K\subset \{1,\ldots, 8\}\setminus \{i\},\;|J|=3,\;|K|=4,\;J\cap K = \emptyset$ \\  \hline
      & $\displaystyle D=8L - 4E_i - 2 \sum_{j\in J} E_j - 3 \sum_{k\in K} E_k,\;i\in \{1,\ldots,8\},\; J,K\subset \{1,\ldots, 8\}\setminus \{i\},\;|J|=3,\;|K|=4,\;J\cap K = \emptyset$ \\  \hline
      & $\displaystyle D=8L - E_i - 3 \sum_{j\in J} E_j,\;i\in \{1,\ldots,8\},\;J=\{1,\ldots,8\}\setminus \{i\} $ \\  \hline
      & $\displaystyle D=9L - 2 E_i - 4 \sum_{j\in J} E_j - 3 \sum_{k\in K} E_k,\; i\in \{1,\ldots,8\},\; J,K\subset \{1,\ldots, 8\}\setminus \{i\},\;|J|=2,\;|K|=5,\;J\cap K = \emptyset$ \\  \hline
      & $\displaystyle D=10L - 3 \sum_{i\in I} E_i - 4 \sum_{j\in J},\; I,J\subset \{1,\ldots,8\},\;|J|=|K|=4,\; J\cap K = \emptyset$ \\  \hline
      & $\displaystyle D=11L - 3E_i - 4 \sum_{j\in J} E_j,\; i\in \{1,\ldots,8\},\; J = \{1,\ldots,8\}\setminus \{i\}$ \\  \hline
\end{tabular}%
}\caption{}
\label{tab_canonical_bundles}

\end{table}

 \end{proposition}

\begin{proof} This is a classical fact that can be found in \cite[\S 2]{Der06} (see also \cite[Table 2]{TVAV09}). For the reader's convenience, we give a self-contained proof.

  Lemma \ref{Numerical characterization of conic bundle} allows us to classify for each $Z_d$, with $1\leq d\leq 8$, the canonical classes $[D]\in \NS(Z_d)$ that produces conic bundles.
Indeed, since $Z_d$ is the blow-up of $\mathbf{P}^2$ in $9-d$ points in general position, 
then \[\operatorname{Pic}(Z_d)=\mathbf{Z}L\oplus \bigoplus_{i=1}^{9-d}\mathbf{Z}E_i\]  where $L$ is the class of pull-back of a line and each $E_i$ is an exceptional divisor. Thus, $K_{Z_d}=-3L+E_1+\dots+E_{9-d}$ and $D=\ell L+a_1E_1+\dots+a_{9-d}E_{9-d}$, so the numerical conditions that impose Lemma \ref{Numerical characterization of conic bundle} are 
$$\begin{cases}
    \ell^2=a_1^2+\dots+ a_{9-d}^2\\ -3\ell+2=a_1+\dots+a_{9-d}
\end{cases}.$$
For each $d$, one gets the possibilities for $D$ of Table \ref{tab_canonical_bundles} solving the Diophantine equation system. A straightforward combinatorial computation gives the total number of conic bundles $N$ for each $d$ depending on $n$.
\end{proof}

\begin{remark}
    If $Z\simeq \mathbf{P}^1\times \mathbf{P}^1$ the only conic bundles are the 2 different projections onto $\mathbf{P}^1$.
\end{remark}

In what follows, we will relate the presence of conic bundles and elliptic fibration, therefore we recall the definition of the latter.
\begin{definition}\label{def:elliptic fibration}
    An elliptic fibration on a smooth projective surface $Z$ is a surjective morphism $\pi:X\rightarrow B$, where $B$ is a smooth algebraic curve, the fiber $\pi^{-1}(t)$ is a curve of genus 1 for all but finitely many $t\in B$ and $\pi$ is relatively minimal.
\end{definition}

\begin{remark}
    Observe that in Definition \ref{def:elliptic fibration} we do not ask for the presence of a section, contrary to other works, e.g. \cite{GarSal20}. 
\end{remark}
\subsection{Non-symplectic involutions on K3 surfaces}

We refer the reader to \cite{HuybrechtsK3,Kondo} for the notation and preliminaries on K3 surfaces and elliptic fibrations (e.g. we consider the ADE lattices to be definite negative).

Let $X$ be a projective K3 surface and $\varphi:X\to X$ be an automorphism of finite order $n$. Let us recall that $\varphi$ is symplectic (resp. purely non-symplectic) if $\varphi^\ast \omega = \omega$ (resp. $\varphi^\ast \omega = \lambda \omega$ for some $\lambda$ a primitive $n$-rooth of the unity), where $\omega$ is a non-degenerate holomorphic $2$-form on $X$ (i.e., $\operatorname{H}^0(X,\Omega_X^2)=\mathbf{C} \cdot \omega$).
Given a K3 surface $X$ and a non-symplectic involution $\iota$ on $X$, we denote $X^\iota$ the fixed locus of the involution and $\NS(X)^{\iota^*}$ the invariant lattice, i.e.
\[\NS(X)^{\iota^*}=\{x\in\NS(X): \iota^*x=x\}.\]

Classical results on non-symplectic involutions on K3 surfaces allow us to classify their fixed loci and invariant lattices. It is an elementary but very useful observation that $H^2(X,\mathbf{Z})^{\iota^*}\subset \NS(X)$ as lattices for any non-symplectic involution $\iota$.

\begin{proposition}{\cite[Theorem 4.2.2]{Nik83}}\label{prop: fixed locus of non-sym inv}
Let $X$ be a K3 surface and $\iota$ be a non-symplectic involution with a non-empty fixed locus. The fixed locus of $\iota$ can be either the disjoint union of two elliptic curves, $X^\iota = E_1 \sqcup E_2$, or the disjoint union $$X^\iota=C_g \sqcup R_1\sqcup \ldots \sqcup R_k,$$ where $C_g$ is a smooth curve of genus $g$ and the $R_i$'s are rational curves. Moreover, $g=\frac{22-r-a}{2}$ and $k=\frac{r-a}{2}$, where $r=\rk \NS(X)^{\iota^*}$ and $a$ is the length of $ \NS(X)^{\iota^*}$, i.e., $2^a= |\discr \NS(X)^{\iota^*}|$.
\end{proposition}

\begin{remark}\label{remark:(r,a)=(10,10)} In the case $X^\iota = \emptyset$ it follows from \cite[Theorem 4.2.2]{Nik83} that necessarily $(r,a,\delta)=(10,10,0)$. Moreover, in that case \cite[Theorem 4.2.4, Proposition 4.2.5]{Nik83} imply that $\NS(X)^{\iota^\ast}\simeq U(2)\oplus E_8(2)$ and the group $\operatorname{Aut}(X)$ is infinite. Since $X^\iota = \emptyset$, we have that $X\slash \langle \iota \rangle$ is an Enriques surface.
\end{remark}

\begin{convention}\label{convention:generic}
    In what follows, we will assume that $X$ is \textit{generic} among the K3 surfaces admitting a non–symplectic involution with a given fixed locus. This is equivalent to the condition the action of $\iota^*$ is trivial on $\NS(X)$. 
\end{convention}

As observed in \cite{Nik83}, the invariants $(r,a)$ allow to recover the topological invariants $(g,k)$ of the fixed locus $X^\iota$. Viceversa,  if the pair $(g,k)$ is known,
a third invariant $\delta$ is needed to identify uniquely $\NS(X)^{\iota^*}$: 
$\delta=0$ if and only if $X^{\iota}\sim 0 \text{ mod } 2$ in $H_2(X,\mathbf{Z})$, otherwise $\delta=1$.
The picture of all possible invariant $(r,a,\delta)$ is presented in \cite[Figure 1]{AleNik89}. For completeness, we include it in Figure \ref{FM}.

\begin{figure}[h]
\centering
\begin{tikzpicture}[scale=.43]
\filldraw [black] 
(1,1) circle (2.5pt) 
(2,2) circle (2.5pt) 
(3,1) circle (2.5pt)
(3,3) circle (2.5pt)
(4,2) circle (2.5pt)
(4,4) circle (2.5pt)
(5,3) circle (2.5pt)
(5,5) circle (2.5pt)
(6,4) circle (2.5pt)
(6,2) circle (2.5pt)
(6,6) circle (2.5pt)
(7,3) circle (2.5pt)
(7,5) circle (2.5pt)
(7,7) circle (2.5pt)
(8,2) circle (2.5pt)
(8,4) circle (2.5pt)
(8,6) circle (2.5pt)
(8,8) circle (2.5pt)
(9,1) circle (2.5pt)
(9,3) circle (2.5pt)
(9,5) circle (2.5pt)
(9,7) circle (2.5pt)
(9,9) circle (2.5pt)
(10,2) circle (2.5pt)
(10,4) circle (2.5pt)
(10,6) circle (2.5pt)
(10,8) circle (2.5pt)
(10,10) circle (2.5pt)
(11,1) circle (2.5pt)
(11,3) circle (2.5pt)
(11,5) circle (2.5pt)
(11,7) circle (2.5pt)
(11,9) circle (2.5pt)
(11,11) circle (2.5pt)
(12,2) circle (2.5pt)
(12,4) circle (2.5pt)
(12,6) circle (2.5pt)
(12,8) circle (2.5pt)
(12,10) circle (2.5pt)
(13,3) circle (2.5pt)
(13,5) circle (2.5pt)
(13,7) circle (2.5pt)
(13,9) circle (2.5pt)
(14,4) circle (2.5pt)
(14,6) circle (2.5pt)
(14,8) circle (2.5pt)
(15,3) circle (2.5pt)
(15,5) circle (2.5pt)
(15,7) circle (2.5pt)
(16,2) circle (2.5pt)
(16,4) circle (2.5pt)
(16,6) circle (2.5pt)
(17,1) circle (2.5pt)
(17,3) circle (2.5pt)
(17,5) circle (2.5pt)
(18,2) circle (2.5pt)
(18,4) circle (2.5pt)
(19,1) circle (2.5pt)
(19,3) circle (2.5pt)
(20,2) circle (2.5pt)
;
\node at (2,0) {$\square$};
\node at (6,2) {$\square$};
\node at (10,0) {$\square$};
\node at (14,2) {$\square$};
\node at (18,0) {$\square$};
\node at (2,2) {$\square$};
\node at (10,2) {$\square$};
\node at (18,2) {$\square$};
\node at (6,4) {$\square$};
\node at (10,4) {$\square$};
\node at (14,4) {$\square$};
\node at (10,6) {$\square$};
\node at (10,8) {$\square$};
\node at (10,10) {$\square$};
\node at (14,6) {$\square$};
\node at (18,4) {$\square$};

\draw[->] (0,0) -- coordinate (x axis mid) (23,0);
    \draw[->] (0,0) -- coordinate (y axis mid)(0,13);
    \foreach \x in {0,1,2,3,4,5,6,7,8,9,10,11,12,13,14,15,16,17,18,19,20}
        \draw [xshift=0cm](\x cm,0pt) -- (\x cm,-3pt)
         node[anchor=north] {$\x$};
          \foreach \y in {1,2,3,4,5,6,7,8,9,10,11}
        \draw (1pt,\y cm) -- (-3pt,\y cm) node[anchor=east] {$\y$};
    \node at (23, -1) {$r$};
    \node at (-1,13) {$a$};

\node at (18,11) {$\bullet$};
\node at (18,10) {$\square$};
\node at (20,11) {$\delta=1$};
\node at (20,10) {$\delta=0$};

   \end{tikzpicture} 
\caption{All possible invariants $(r,a,\delta)$.}
  \label{FM}
\end{figure}

\vspace{2mm}

The following result characterizes the quotient of a K3 surfaces by a non-symplectic automorphism.

\begin{proposition}\label{Prop: quotient of K3 is rational or Enriques}
Let $X$ be a projective K3 surface and $\varphi$ be a non-symplectic automorphism of finite order. Then, \begin{enumerate}[i)]
    \item the quotient $X/\langle \varphi \rangle$ is a rational surface or birational to an Enriques surface;
    \item if $\varphi$ is an involution, then the fixed locus of $\varphi$ is empty if and only if $X/\langle \varphi \rangle$ is an Enriques surface.
\end{enumerate}

\end{proposition}

\begin{proof} See \cite[Lemma 4.8]{HuybrechtsK3} for $i)$. For $ii)$, the ``only if'' part is proven in \cite[Lemma 1.2]{Zha98}. 
Suppose that $\varphi$ is an involution. We know that if the fixed locus is non-empty, then $X^\varphi=D$ has codimension one and it is a disjoint union of smooth curves. Furthermore, the canonical divisor of cyclic coverings implies that $D \sim 2K_{X/\langle \varphi \rangle}$ and so the quotient is not an Enriques surface. Suppose that the fixed locus $X^\varphi$ is empty. Then, $2 \pi^* K_{X/\langle \varphi \rangle} \sim K_X \sim 0$, implying that $X/\langle \varphi \rangle$ is an Enriques surface.
\end{proof}

\begin{remark}
    It is a classical fact that K3 surfaces admitting a non-symplectic automorphism of finite order are projective and the N\'eron-Severi lattice is hyperbolic. See \cite[Chapter 15, Corollary 1.10]{HuybrechtsK3} and \cite[Theorem 3.1]{Nik79} for instance.
\end{remark} 

\section{Strictly elliptic involutions}
Let $X$ be a K3 surface and $\iota$ a non-symplectic involution on $X$.
According to the notation in \cite[\S 2.8]{AleNik06}, it is natural to classify $\iota$ into three categories: $\iota$ is of \textit{elliptic} type if $X^\iota$ 
contains a curve of genus $g\geq 2$. This is equivalent to $r+a\leq18$, $(r,a,\delta)\neq (10,8,0)$. The involution
$\iota$ is of \textit{parabolic} type if $X^\iota$ contains a genus 1 curve, which is equivalent to $r+a=20$ or $(r,a,\delta)= (10,8,0)$. Finally, $\iota$ is of \textit{hyperbolic} type if $X^\iota=\emptyset$ or $X^\iota$ only contains rational curves, which is equivalent to $r+a=22$ or $(r,a,\delta)= (10,10,0)$. 

In this work, we specifically focus on the case of elliptic type, and in particular, we consider K3 surfaces that admit a non-symplectic involution with no rational curves in the fixed locus.

\begin{definition}\label{def: strictly elliptic type} Let $X$ be a K3 surface and $\iota$ be a non-symplectic involution on $X$. We say that $\iota$ is of strictly elliptic type if its fixed locus is given by $X^\iota=C_g$, where $C_g$ is a smooth irreducible curve of genus $g\geq 2$.
\end{definition}

\begin{remark}\label{NS} 
Looking at Figure \ref{FM}, one can observe that the N\'eron-Severi groups of the strictly elliptic type K3 surfaces correspond to points on the line $r=a$ with $r\leq 9$ and they all have $\delta=1$ except for the case $r=2$, where both $\delta=0$ and $\delta=1$ are possible.
As a consequence of \cite{Nik83}, possibilities for $\NS(X)$ are as follows: if $r=2, \delta=0$, then $\NS(X)\iso U(2)$; otherwise if $\delta=1$ one has 
$$\NS(X) \iso \begin{cases}
\langle 2\rangle, \text{ if } r=1\\
\langle 2\rangle\oplus A_1, \text{ if } r=2 \\
U(2)\oplus A_1^{\oplus r -2}, \text{ if } r \geq 3.
\end{cases}$$

\end{remark}

Observe that by \cite[\S 2.8]{AleNik06}, the automorphism group of a K3 surface with a non-symplectic involution of elliptic type is finite.

\begin{remark}\label{remark:remaining-cases}

As mentioned in the Introduction, our methods rely on the correspondence between strictly elliptic K3 surfaces and smooth del Pezzo surfaces, as established in \cite{AleNik06}. This correspondence allows us to analyze the case where the fixed locus corresponds to a single smooth irreducible curve $C_g$ of genus $g\geq 2$. Additionally, prior works \cite{Ogu89,Klo06,CG,GarSal19, GarSal20} are devoted to describe all (jacobian) elliptic fibrations in Figure \ref{FM} except for the cases $(r,a,\delta)=(10,10,0)$ and $(r,a,\delta)=(10,10,1)$. The former case corresponds precisely to the situation $X^\iota = \emptyset$ (see Remark \ref{remark:(r,a)=(10,10)}) and thus $X\slash \langle \iota \rangle$ is an Enriques surface, while the former case arises when the fixed locus is given by a single smooth elliptic curve $E$.

It is noteworthy that the case $(r,a,\delta)=(10,10,1)$ falls beyond our analysis and it should be noted that the nature of the corresponding elliptic fibrations is necessarily different from the strictly elliptic case. More precisely, in \cite[\S 2]{GarSal20} the authors consider non-symplectic involutions $\iota:X\to X$ satisfying the Convention \ref{convention:generic} and classify the possible elliptic fibrations $\mathcal{E}:X\to \mathbf{P}^1$ into two types:
\begin{itemize}
    \item {\bf Type 1.} $\iota$ maps each fiber of $\mathcal{E}$ to itself. 
    \item {\bf Type 2.} $\iota$ maps at least one fiber of $\mathcal{E}$ to another fiber of $\mathcal{E}$. 
\end{itemize}
It follows from the proof of \cite[Proposition 2.5]{GarSal20} that strictly elliptic K3 surfaces only admit elliptic fibrations of Type 1.

On the other hand, in the case $(r,a,\delta)=(10,10,1)$ we have that the quotient $Z=X\slash \langle \iota \rangle$ is a smooth rational surface with $K_Z^2=0$ and $-2K_Z \sim E$, where $X^\iota \simeq E$ is an elliptic curve (cf. proof of Proposition \ref{prop: K3 of strictly elliptic type} below). In particular, it follows from \cite[Proposition 2.2]{CD12} that there exist an irreducible pencil of sextics in $\mathbf{P}^2$ with 9 nodes $p_1,\ldots,p_9$ as base points such that $Z\cong \operatorname{Bl}_{p_1,\ldots,p_9}(\mathbf{P}^2)$, and the elliptic fibration $\phi_{|E|}:Z\to \mathbf{P}^1$ is the proper transform of this pencil.

Finally, observe that by \cite[Corollary 3.3]{CliMal23}, the case $(r,a,\delta)=(6,4,0)$ does not admit jacobian elliptic fibrations and therefore it is not considered in the analysis of \cite{GarSal20} despite the presence of a rational curve in $X^\iota$. On the other hand, by \cite[Table 9]{Nak07}, $Z=\operatorname{Proj}\oplus_{m\geq 0}\operatorname{H}^0(X,\mathcal{O}_X(mC_g))^\iota$ is a \emph{singular} del Pezzo surface of Gorenstein index two that can be realized via a Sarkisov link starting from $\mathbf{P}(1,1,4)$ or as a hypersurface of degree $5$ in $\mathbf{P}(1,1,1,4)$ (see \cite[Proposition 7.4, 7.11]{Nak07}).
\end{remark}

\begin{proposition}\label{prop: K3 of strictly elliptic type}
   Let $X$ be a K3 surface admitting an involution of strictly elliptic type. Suppose that $C_g$ is the smooth irreducible curve of genus $2\leq g\leq 10$ fixed by $\iota$.
   \begin{enumerate}
       \item  If $\delta=0$, then the K3 surface $X$ can be realized as the double cover of $\mathbf{P}^1\times \mathbf{P}^1$ over a smooth irreducible curve $C\in |\mathcal{O}_{\mathbf{P}^1 \times \mathbf{P}^1}(4,4)|$ of genus $9$.
       \item If $\delta=1$, then the K3 surface $X$ can be realized as the double cover of a del Pezzo surface $Z_{d}$ of degree $d = g - 1$ over a smooth irreducible curve $C$ belonging to the linear system $|-2K_{Z_d}|$.
   \end{enumerate}

  Moreover, in the case $\delta=1$ the image of $C$ in $\mathbf{P}^2$ is a nodal sextic curve $\Gamma_d$ that passes through $9 - d$ points corresponding to the points where the blow-up of $\mathbf{P}^2$ is done to obtain $Z_d$. 
\end{proposition}

\begin{proof}
    By Proposition \ref{Prop: quotient of K3 is rational or Enriques} and the fact that $X^\iota$ is of pure codimension 1 (see Proposition \ref{prop: fixed locus of non-sym inv}), the quotient $X/\langle \iota \rangle$ is a smooth rational surface. Furthermore, since the fixed locus of $\iota$ is an irreducible curve $C_g$ of genus $g\geq 2$, the quotient $X/\langle \iota \rangle$ is a del Pezzo surface of degree $d=K^2_{X/\langle \iota \rangle}=g-1$. Denote by $Z=X/\langle \iota \rangle$ and by $C$ the image under the quotient map of the curve $C_g$. We have that $C\in |-2K_{Z}|$, and by the genus formula, $g(C)=g(C_g)=K_{Z}^2+1=d+1$. According to the classification of strictly involutions of elliptic type, we know that $g$ is an integer number at most $10$, hence the degree of $Z_d$ satisfies $1\leq d \leq 9$. Consequently, either $Z\simeq \mathbf{P}^1 \times \mathbf{P}^1$ (and $d=8$) or $Z\simeq Z_d$ is obtained as the blow-up of $\mathbf{P}^2$ at some generic points $p_1,\dots,p_{9-d}$. In the latter case, let $\beta$ be the blow-up map and $E_i's$ be the corresponding exceptional divisors, then  $C\in |-2K_{Z_d}|=|6L-2\sum_{i=1}^{9-d}E_i|$, where $L$ is the pull-back of a line in $\mathbf{P}^2$ via $\beta$. This implies that $C$ can be considered as the strict transform of an irreducible curve $\Gamma_d \subseteq \mathbf{P}^2$ of degree $6$ passing with multiplicity two at each point $p_1,\dots,p_{9-d}$. It is worth noting that curves with these properties can always be found; see e.g. \cite[Chapter 5, Section 5.2, Theorem 1]{Fulton}.

    Finally, the fact that $\delta=0$ corresponds precisely to double covers of $\mathbf{P}^1\times \mathbf{P}^1$ follows from Proposition \ref{prop: fixed locus of non-sym inv}, Table \ref{FM}, and the topological condition \cite[Formula (38), page 32]{AleNik06} imposing that $\frac{1}{4}[C]$ is an integral homology class in $H_2(X\slash \langle \iota \rangle,\mathbf{Z})$.
\end{proof}

\begin{remark}
    It is worth noticing that the case where $X^\iota$ consist of a curve of genus 10 fits in the framework of the previous proposition when considering $d=9$, i.e., the rational surface is $Z_9\simeq \mathbf P^2$. In this case, as we will see in Theorem \ref{thm_main}, there are no elliptic fibrations on $X$, therefore our analysis will focus on $d\leq 8$.
\end{remark}

Following the construction presented  in Proposition \ref{prop: K3 of strictly elliptic type}, if we take the double cover of the plane $\mathbf P^2 $ branched at a nodal sextic curve $\Gamma_d$ we get a singular surface $Y$, with ADE singularities, such that its minimal resolution $X$ is a K3 surface, and the strict transform of $\Gamma_d$ is a smooth curve $C_g$ on $X$. So, we have a diagram as below.

\begin{figure}[H]
    \centering
    \begin{tikzcd}
    C_g \subset \arrow[d, "1:1", maps to] & X \arrow[d, "2:1"'] \arrow[r, "\text{min. res}"] & Y \arrow[d, "2:1"] &               \\
    C \subset                             & Z_d \arrow[r, "\beta"]                           & \mathbf P^2        & \supset\Gamma_d
    \end{tikzcd}
    \caption{Strictly elliptic K3 surfaces and smooth del Pezzo surfaces as in Proposition \ref{prop: K3 of strictly elliptic type} (2).}
    \label{fig:K3 elliptic type}
\end{figure}

\begin{remark} It is worth mentioning that in \cite[\S 7]{GarSal20} the authors consider a K3 surface $X$ together with a non-symplectic involution $\iota:X\to X$ as in Convention \ref{convention:generic} and such that $X^{\iota}=C_g \sqcup R_1\sqcup \ldots \sqcup R_k$ consists of a smooth irreducible curve $C_g$ of genus $g\geq 2$ and $k\geq 1$ rational curves. In this setting, some of these K3 surfaces are such that $X\slash \langle \iota \rangle \simeq \mathbf{P}^2$, and it is observed that in such instances, $X$ can be realized as the double covering of $\mathbf{P}^2$ along a \emph{reducible} sextic curve. We refer the reader to the Table in \cite[\S 7.1]{GarSal20} for further details (see also \cite[\S 5]{GarSal19} and \cite[\S 5.10]{CG}).
\end{remark}

\begin{lemma}\label{lem: Fibras conic bundle en Zd}
Let $Z$ be a del Pezzo surface of degree $d$ and $f\colon Z\to \mathbf{P}^1$ be a conic bundle on $Z$ with $F_{x}$ the fiber over a point $x\in\mathbf{P}^1$. If $E$ is an irreducible curve such that $E\subseteq \operatorname{Supp}(F_x)$, then $E^2\in \{-1,0\}$ and $g(E)=0$. Let\footnote{The existence of such a curve can be deduced from the Smooth Divisor Theorem (see \cite[Theorem 1.5]{AleNik06}), or simply from the classification of smooth del Pezzo surfaces. Indeed, the divisor $-2K_{Z_d}$ is very ample for $d\geq 2$ and it defines a double cover $\phi_{|-2K_Z|}:Z\to Q$ for $d=1$, where $Q\subseteq \mathbf{P}^3$ is a quadric cone. See \cite[Chapter IV]{Bea96} for details.} $C\in |-2K_{Z}|$ be a smooth irreducible curve, then 
    \begin{equation*}
E^2=
    \begin{cases}
        -1 &\text{ if } E\cdot C=2,\\
        \,\,\,\,\, 0 & \text{ if } E \cdot C=4.
    \end{cases}
\end{equation*} 
\end{lemma}

\begin{proof}

By the correspondence in Lemma \ref{Numerical characterization of conic bundle}, we can associate to any conic bundle $f:Z\to\mathbf{P}^1$ a conic class, i.e., an effective class $[D]\in \operatorname{NS}(Z)$ such that $D^2=0, K_{Z}\cdot D=-2$, and $f=\phi_{|D|}$ is corresponding the induced map. Furthermore, if $E$ is an irreducible component of $\operatorname{Supp}(F_x)$, we have that $E\cdot F_x=E\cdot D=0$. Suppose that $E^2>0$. By the Hodge index theorem, $D$ is numerically trivial. However, this contradicts $K_{Z}\cdot D=-2$. Thus, $E^2\leq0$.

By the genus formula, we have that $2p_a(E)-2=E^2+K_{Z}\cdot E \leq 0$ and thus $2p_a(E)-2<0$, since $E^2\leq 0$ and $-K_Z\cdot E>0$ by Nakai-Moishezon ampleness criterion. Hence, $p_a(E) = 0$, and thus $E \simeq \mathbf{P}^1$. Finally, for a smooth irreducible curve $C\subseteq Z$ such that $-2K_{Z}\sim C$ we have that $E \cdot C = 2E^2 + 4$, from which the last statement of the lemma follows thoroughly.
\end{proof}

\begin{remark}\label{rem_singular_fibers}
Note that the smooth fibers of any such conic bundle over del Pezzo surfaces are rational curves intersecting the curve $C\in |-2K_{Z}|$ at four points, taking multiplicities into account. Additionally, the singular fibers consist of two rational curves with self-intersection $-1$ that intersect at a single point and each component intersects $C$ at two points.  
\end{remark}

We can now state our main theorem concerning the correspondence between conic bundles on $Z$ and elliptic fibrations on $X$.

\begin{theorem} \label{thm_main}
Let $(X,\iota)$ be a pair of strictly elliptic type and ${Z=X\slash \langle \iota \rangle}$ be the quotient smooth del Pezzo surface as in Proposition \ref{prop: K3 of strictly elliptic type}, and let $\pi:X\to Z$ be the quotient map. 

If $d=9, g=10$, the K3 surface $X$ does not admit any elliptic fibration. If $d\leq 8$, there is a correspondence
\[
\left\{\begin{array}{c}
     \textup{Conic bundles}\\ f:Z\to \mathbf{P}^1
\end{array} \right\} \xrightarrow{\;\;\sim \;\;} \left\{ \begin{array}{c}
     \textup{Elliptic fibrations}\\ \mathcal{E}:X\to \mathbf{P}^1
\end{array} \right\},\;f\longmapsto f\circ \pi
\]
Moreover, the fixed curve $X^\iota=C$ is a bisection of $\mathcal{E}$, i.e., $E\cdot C=2$ for the general fiber $E$ of $\mathcal{E}:X \rightarrow \mathbf{P}^1$.
\end{theorem}

\begin{proof}

    The case $Z\simeq \mathbf{P}^1 \times \mathbf{P}^1$ is treated in \cite[Proposition 2.4]{Dol73} and \cite[\S 3]{Reid}, so we will assume that $Z\simeq Z_d$ is the blow-up of $9-d$ points in $\mathbf{P}^2$ in general position (cf. \cite[Theorem 3.2]{Reid} where the case of $Z_8 \simeq \mathbf{F}_1$ is also considered).

We first observe that when $g=10$, the K3 surface $X$ does not admit elliptic fibrations. Indeed  as observed in Remark \ref{NS}, in this case $\NS(X)\iso\langle 2\rangle$, i.e., it is generated by a class of square 2. Thus there are no classes $D\not\sim 0$ with $D^2=0$ and therefore no elliptic fibrations.

    First, let us consider $\mathcal{E}:X\to \mathbf{P}^1$ an arbitrary elliptic fibration. Since $\iota$ is strictly elliptic, we know that $\iota$ must map each fiber of $\mathcal{E}$ to itself (see Remark \ref{remark:remaining-cases}), and thus $\mathcal{E}$ factors through a fibration $f:X\slash \langle \iota \rangle \simeq Z_d \to \mathbf{P}^1$. Since $Z_d$ is a smooth del Pezzo surface and $f$ has connected fibers, it follows that $f$ is a $K_{Z_d}$-negative contraction with one-dimensional fibers and thus is a conic bundle by \cite[ Theorem 3.1 (ii)]{Ando}.

    Conversely, if we consider a conic bundle $f: Z_d \to \mathbf{P}^1$ and define $\mathcal{E}:=f\circ \pi:X\to \mathbf{P}^1$. Here, we have that $\pi(C)\in |-2K_{Z_d}|$. In particular, if we denote by $F$ the general fiber of $f:Z_d\to \mathbf{P}^1$ we have that $\pi(C)\cdot F=4$ and then it follows from the Riemann-Hurwitz formula that $\mathcal{E}:X\to \mathbf{P}^1$ is a elliptic fibration. The induced conic bundle is precisely $f:Z_d \to \mathbf{P}^1$ and thus we get the desired correspondence. Finally, the fact that $C$ is a bisection of $\mathcal{E}:X\to \mathbf{P}^1$ follows directly from the projection formula. 
\end{proof}

\begin{figure}[H]
    \centering
    \begin{tikzcd}
    X \arrow[r, "\Ee"] \arrow[d, "\pi"'] & \mathbf{P}^1 \\
    Z \arrow[ru, "\exists f"']           &             
    \end{tikzcd}
    \caption{Elliptic fibrations on strictly elliptic K3 surfaces and conic bundles on del Pezzo surfaces.}
    \label{fig: elliptic fibration vs conic bundle}
    \end{figure}

\begin{example}\label{example:bisections}
    Let $Z_d$ be a del Pezzo surface of degree $d\leq 8$ obtained as the blow-up of $9-d$ points in general position in $\mathbf{P}^2$. Fix a conic bundle $f:Z_d\to \mathbf{P}^1$ and let $\mathcal{E}:X\to \mathbf{P}^1$ be the corresponding elliptic fibration.

    Let $F$ be the general fiber of $\mathcal{E}:X\to \mathbf{P}^1$, $\Gamma$ be the general fiber of $f:Z_d\to \mathbf{P}^1$ and $E\subseteq Z_d$ be a $(-1)$-curve. By the projection formula 
    \[
    F\cdot \pi^\ast(E)=\deg(\pi) (\Gamma \cdot E) = \deg(\pi)=2 \textup{ as long as }\Gamma\cdot E = 1,
    \]
  and hence these $(-1)$-curves on $Z_d$ induce bisections of $\mathcal{E}:X\to \mathbf{P}^1$. The condition $\Gamma \cdot E = 1$ can be explicitly verified by means of Theorem \ref{thm:(-1)-curves} and Proposition \ref{Classes of conic bundles on Z}. For instance, following Convention \ref{convention:Z_d}, in $Z_4$ we can consider $E=2L-E_1-E_2-E_3-E_4-E_5$ and we can check that:
  \begin{itemize}
      \item If $\Gamma = L - E_1$ then $E\cdot \Gamma = 1$.
      \item If $\Gamma = 2L - E_1 - E_2 - E_3 - E_4$ then $E\cdot \Gamma = 0$.
  \end{itemize}
\end{example}

\section{Examples and Néron-Severi lattices}

We will now present some examples.

\begin{example}
    Let $X$ be the smooth quartic surface  given by the equation $$\{x_0^4+x_1^4+x_2^4+x_3^4+6(x_0^2x_1^2-x_1^2x_2^2+x_1^2x_3^2+x_0^2x_2^2-x_0^2x_3^2+x_2^2x_3^2)=0\}\subseteq\mathbf{P}^3.$$
    We consider the automorphism
    \[
    \iota:X\to X,\;[x_0:x_1:x_2:x_3]\mapsto [x_0:x_1:x_2:-x_3]
    \]
    which is a non-symplectic involution on $X$, and the fixed locus of $\iota$ is the smooth quartic curve
$C=\{x_0^4+x_1^4+x_2^4+6(x_0^2x_1^2-x_1^2x_2^2+x_0^2x_2^2)=0\}\subseteq \mathbf{P}^2$ of genus $g(C)=3$. 

The quotient $X/\langle \iota \rangle=Z$ is a smooth rational surface, and  by the formula of the canonical divisor of double coverings (see \cite[Chapter V, \S 22]{BHPV04}) we have that $C\sim-2K_Z$. Therefore, the bi-anticanonical divisor $-2K_Z$ is ample and $K^2_Z=2$, i.e., $Z\simeq Z_2$ a del Pezzo surface of degree two. Conversely, a del Pezzo surface of degree $2$ is isomorphic to a double cover of $\mathbf{P}^2$ branched at a smooth quartic curve (see e.g. \cite[Chapter IV]{Bea96}).

Note that \cite[\S 5.3]{CD24} gives a correspondence between plane sextic curves with $7$ nodes in general position and smooth planar quartic curves. In particular, the smooth quartic curve $C\subseteq \mathbf{P}^2$ corresponds a sextic curve $\Gamma\subseteq \mathbf{P}^2$ with $7$ nodes in general position. By considering the blow-up of $\mathbf{P}^2$ at the $7$ nodes of $\Gamma$, followed by taking the double cover ramified at the strict transform of $\Gamma$, we obtain the K3 surface $X$. 

The surface $X$ has N\'eron–Severi lattice $U(2)\oplus A_1^{\oplus 6}$, as outlined in \cite[Sect. 8.4]{Rou}. It is noteworthy, in accordance with \cite[Theo. 8.2]{avila2023invariant}, that $X$ corresponds to Burnside's conjecture as the most symmetric smooth quartic surface. Specifically, its group of projective automorphisms is $\mathbf{Z}_2^4\cdot \mathfrak{S}_5$.
\end{example}

As we observed, K3 surfaces of strictly elliptic type have finite automorphism group (see \cite[\S 2.8]{AleNik06}) and thus they are related with the works of Roulleau, Artebani and Correa \cite{Rou, ACR}.

For instance, if $d=7$, the corresponding K3 surface has Picard rank 3 and is obtained as double cover of $\mathbf P^2$ branched over a sextic curve with two nodes $p_1, p_2$. By \cite[Prop. 3.4]{Rou} the surface admits three $(-2)$-curves. 
Two of them are contracted to $p_1$ and $p_2$, while the image of the third is the line through the two nodal points.

Similar descriptions and properties regarding the configuration of the $(-2)$-curves for $1\leq d\leq 6$ are provided in \cite{Rou, ACR}. We have summarized the references for each case in Table \ref{ex-roulleau}. 

\begin{table}[h]
\begin{tabular}{|l|l||l|}
\hline
$r$&$d$ & 
    $\operatorname{NS}(X)$
     \\ \hline
$3$&$7$ & $S_{1,1,2}\simeq U(2)\oplus A_1$,  \cite[Sect. 3.4]{Rou}\\ \hline
$4$&$6$ &  $U(2)\oplus A_1^{\oplus 2}$,  \cite[Prop. 2.11]{ACR}             \\ \hline
$5$&  $5$  & $U(2)\oplus A_1^{\oplus 3}$, \cite[Sect. 5.2]{Rou}               \\ \hline
$6$& $4$   &  $U(2)\oplus A_1^{\oplus 4}$, \cite[Sect. 6.7]{Rou}              \\ \hline
$7$& $3$   &    $U(2)\oplus A_1^{\oplus 5}$,  \cite[Sect. 7.3]{Rou}            \\ \hline
$8$&  $2$  &   $U(2)\oplus A_1^{\oplus 6}$, \cite[Sect. 8.4]{Rou}             \\ \hline
$9$&   $1$  &   $U(2)\oplus A_1^{\oplus 7}$, \cite[Sect. 9.5]{Rou}             \\ \hline
\end{tabular}
\caption{}
\label{ex-roulleau}
\end{table}

In the case $d=8$ (and subsequently $r=2$), we distinguish between the cases $\delta=0$ and $\delta=1$, where $Z\simeq \mathbf{P}^1\times \mathbf{P}^1$ and $Z \simeq \mathbf{F}_1$ respectively. In the $\delta=0$ case, the Néron–Severi lattice $\NS(X)$ is generated by the pullback via $\pi:X\to Z$ of the two fibers of the canonical projections, and thus $\NS(X)\simeq U(2)$. In the $\delta=1$ case, $\NS(X)$ is generated by the pullback via $\pi:X\to Z$ of the unique $(-1)$-curve in $\mathbf{F}_1$ and the class in $\mathbf{F}_1$ obtained by the pullback of a line in $\mathbf{P}^2$ via the contraction $\mathbf{F}_1\to \mathbf{P}^2$, yielding $\NS(X)\simeq \langle 2 \rangle \oplus A_1$ (cf. Remark \ref{NS}).

\begin{example}\label{example:P1xP1}

    Let $X$ be the double cover of $\mathbf{P}^1\times \mathbf{P}^1$ branched at a smooth bi-quartic curve $C=\{f_{4,4}([x,y],[s,t])=0\}$, i.e.,
    \[X:=\{w^2=f_{4,4}([x,y],[s,t])\}.\]

We will exhibit an elliptic fibration $\mathcal{E}:X\rightarrow\mathbf{P}^1$
 that admits fibers of types $I_0,I_1,I_2$, $II$ and $III$. 
Let $C$ be the smooth bi-quartic curve in $\mathbf{P}^1\times \mathbf{P}^1$ given by the equation
\begin{equation*}\begin{split}
    C=\{&x^4s^4-\tfrac{1}{2}x^3ys^3t+\tfrac{1}{2}x^3ys^2t^2-x^2y^2s^4-\tfrac{1}{2}x^2y^2s^3t+\tfrac{3}{2}x^2y^2s^2t^2\\&-x^2y^2t^4+\tfrac{1}{2}xy^3s^2t^2+y^4z^4-2y^4s^2t^2+y^4t^4=0\}.
\end{split}\end{equation*}

We consider the elliptic fibration $\mathcal E:X\rightarrow\mathbf P^1$ obtained composing with 
the first projection onto $\mathbf P^1$. In the affine chart $t=1$, the fiber over a generic point $[a:b]\in \mathbf{P}^1$ is given by 
    \begin{equation*}
        w^2=(a^4-a^2b^2+b^4)s^4+\tfrac{1}{2}(-a^3b+a^2b^2)s^3+\tfrac{1}{2}(a^3b+3a^2b^2+ab^3-4b^4)s^2+(\pc{-}a^2b^2+b^4).
\end{equation*}    
Using classical  Weierstrass' methods (see e.g. \cite[Chapter 10, Theorem 2]{mordell1969diophantine}) we can see that the generic fiber can be written as $w^2=f_3(z,a,b)$, where $f_3(z)$ is a polynomial of degree three in $z$ such that for generic $[a:b]\in \mathbf{P}^1$ the corresponding curve is irreducible.

If we denote by $L_{[a:b]}$ the line $\{[a:b]\}\times \mathbf{P}^1$, we have the following:
\begin{enumerate}
    \item $C\cap L_{[1:0]}=\{s^4=0\}$, so the fiber over $[1:0]$ is ${w^2=s^4}$, of type $III$.
    \item $C\cap L_{[1:1]}=\{s^4-s^3t=0\}$, so the fiber over $[1:1]$ is ${w^2=s^3(s-t)}$, of type $II$.
    \item $C\cap L_{[0:1]}=\{s^4-2s^2t^2+t^4=0\}$, so the fiber over $[0:1]$ is ${w^2=(s-t)^2(s+t)^2}$, of type $I_2$.
    \item $C\cap L_{[1:-1]}=\{s^4-s^2t^2=0\}$, so the fiber over $[1:-1]$ is ${w^2=s^2(s^2-t^2)}$, of type $I_1$.

\end{enumerate}

As mentioned in the Introduction, the analysis of double covers of $\mathbf{P}^1\times \mathbf{P}^1$ branched over curves of bi-degree $(4,4)$ was undertaken in \cite{Dol73,Reid}, where the authors already observed that the corresponding K3 surface $X$ admits an elliptic fibration (see \cite[Theorem 4.4]{Dol73} and \cite[\S 3]{Reid}). This study builds upon classical works by Enriques and Campedelli.
\end{example}

\section{Singular fibers and bisections}

 Let $(X,\iota)$ and $\pi:X\to Z=X\slash \langle \iota \rangle$ be as in Theorem \ref{thm_main}. 
It is a well-known fact  that if there exists a primitive embedding of the lattice $U\hookrightarrow \NS(X)$, then the K3 surface $X$ admits a jacobian elliptic fibration (see e.g. \cite[Remark \S 11, 1.4]{HuybrechtsK3}). However, according to \cite[Corollary 3.3]{CliMal23}, K3 surfaces of strictly elliptic type do not possess such embeddings, and consequently, they do not admit jacobian elliptic fibrations. On the other hand, since there is an embedding of the lattice $U(2)$ in $\NS(X)$, these surfaces admit bisections. More precisely, one can choose a class $L$ in $U(2)$ with $L^2=0$ and assume that $L$ is nef modulo the action of the Weyl group. By \cite[Remark \S 8, 2.13]{HuybrechtsK3}, $X$ admits an elliptic fibration. Furthermore, since the divisibility of $L$ is either 1 or 2, the elliptic fibration admits either sections or bisections. In our specific context, where elliptic fibrations do not admit sections, we can therefore conclude the existence of bisections.

\begin{proposition}\label{theoretical fibers}

  Let $f:Z\to \mathbf{P}^1$ be a conic bundle on the smooth del Pezzo surface $Z$ and let $\mathcal{E}:X\to\mathbf{P}^1$ be the induced elliptic fibration on $X$. If $C\in|-2K_Z|$ is the branching locus of $\pi:X\to Z$ then: 
\begin{enumerate}

\item If $F$ is a smooth fiber of $f$, then the fibers of $\mathcal{E}$ are of the following types (see Table \ref{table:smooth fibers}): 
\begin{itemize}
    \item $I_0$ if $C$ meets $F$ in 4 distinct points;
    \item $I_1$ if $C$ meets $F$ in 2 simple points and a double point;
    \item $I_2$ if $C$ meets $F$ in two double points;
    \item $II$ if $C$ meets $F$ in two points:  a simple one and a point with multiplicity 3;
    \item $III$ if $C$ meets $F$ in a single point with multiplicity 4.
\end{itemize}
\item If $F=F_1+F_2$ is a singular fiber of $f$, let $P=F_1\cap F_2$. Then the fibers of $\mathcal E$ are of the following types (see Table \ref{table:sing_fibers}):
    \begin{itemize}
        \item $I_2$, if $C$ meets each $F_i$ in two simple points, distinct from $P$;
        \item $I_3$, if $C$ meets $F_2$ in two simple points, and $C$ meets $F_1$ in a double point, distinct from $P$;
        \item $I_4$, if $C$ meets each $F_i$ in a double point, distinct from $P$; 
        \item $III$,  if $P\in C$ and $C$  meets each $F_i=2$ in two simple points
        \item $IV$, if $C$ meets $F_1$ in two simple points, while $C$ meets $F_2$ in $P$ with multiplicity 2.
    \end{itemize}  \end{enumerate}
\begin{table}[h!] 

\setlength{\tabcolsep}{5mm} 
\def\arraystretch{1.25} 
\centering
\makebox[\linewidth]{
\begin{tabular}{|c|c|c|c|c|}
\hline
\hspace{-4mm}
\begin{tikzpicture}[x=0.75pt,y=0.75pt,yscale=-1,xscale=1,scale=0.4]

\draw    (405.55,8.95) -- (407.5,190) ;
\draw [color={rgb, 255:red, 208; green, 2; blue, 27 }  ,draw opacity=1 ]   (380,32) .. controls (412.5,23) and (435.5,28) .. (438.5,64) .. controls (441.5,100) and (382.44,56.59) .. (367.47,94.29) .. controls (352.5,132) and (439.5,103) .. (444.47,131.48) .. controls (449.43,159.95) and (418.5,147) .. (360.5,160) ;
\draw    (285,17) .. controls (292.7,41.59) and (288.05,80.08) .. (271.27,83.54) .. controls (254.5,87) and (260.66,52.85) .. (236.5,52) .. controls (212.34,51.15) and (207.5,140) .. (237.5,138) .. controls (267.5,136) and (260.53,110.33) .. (272.5,112) .. controls (284.47,113.67) and (286.5,160) .. (275.5,181) ;
\draw    (294.5,99) -- (347.55,99) ;
\draw [shift={(349.55,99)}, rotate = 180] [color={rgb, 255:red, 0; green, 0; blue, 0 }  ][line width=0.75]    (10.93,-3.29) .. controls (6.95,-1.4) and (3.31,-0.3) .. (0,0) .. controls (3.31,0.3) and (6.95,1.4) .. (10.93,3.29)   ;

\draw (263,199.4) node [anchor=north west][inner sep=0.75pt]    {$I_{0}$};
\draw (342,197.4) node [anchor=north west][inner sep=0.75pt]  [color={rgb, 255:red, 208; green, 2; blue, 27 }  ,opacity=1 ]  {$[ 1,1,1,1]$};
\draw (315.11,66.4) node [anchor=north west][inner sep=0.75pt]    {$\pi $};

\end{tikzpicture}
 \hspace{-5mm}  & \hspace{-4mm}\begin{tikzpicture}[x=0.75pt,y=0.75pt,yscale=-1,xscale=1,scale=0.4]

\draw    (399.43,1.85) -- (401.03,174.42) ;
\draw [color={rgb, 255:red, 208; green, 2; blue, 27 }  ,draw opacity=1 ]   (345.55,28.27) .. controls (377.53,22.96) and (396.62,20.16) .. (399.82,52.02) .. controls (403.02,83.88) and (374.32,56.35) .. (358.36,89.72) .. controls (342.41,123.09) and (435.14,97.43) .. (440.43,122.63) .. controls (445.73,147.83) and (412.75,136.36) .. (350.93,147.87) ;
\draw    (298.94,97) -- (341.55,97) ;
\draw [shift={(343.55,97)}, rotate = 180] [color={rgb, 255:red, 0; green, 0; blue, 0 }  ][line width=0.75]    (10.93,-3.29) .. controls (6.95,-1.4) and (3.31,-0.3) .. (0,0) .. controls (3.31,0.3) and (6.95,1.4) .. (10.93,3.29)   ;
\draw    (296.91,64.55) .. controls (202.85,189.07) and (195.36,87.29) .. (210.88,66.12) .. controls (226.4,44.95) and (275.12,102.38) .. (295.76,130.75) ;

\draw (309.11,64.4) node [anchor=north west][inner sep=0.75pt]    {$\pi $};
\draw (243,171.4) node [anchor=north west][inner sep=0.75pt]    {$I_{1}$};
\draw (340.27,176.55) node [anchor=north west][inner sep=0.75pt]  [color={rgb, 255:red, 208; green, 2; blue, 27 }  ,opacity=1 ]  {$[ 2,1,1]$};

\end{tikzpicture} \hspace{-5mm}
 & \hspace{-4mm} \begin{tikzpicture}[x=0.75pt,y=0.75pt,yscale=-1,xscale=1,scale=0.4]

\draw    (429,3.25) -- (430.5,154.69) ;
\draw [color={rgb, 255:red, 208; green, 2; blue, 27 }  ,draw opacity=1 ]   (378.45,26.43) .. controls (408.45,21.77) and (426.37,19.32) .. (429.37,47.28) .. controls (432.37,75.24) and (405.44,51.08) .. (390.47,80.36) .. controls (375.5,109.65) and (429,97.42) .. (430,117.61) .. controls (431,137.8) and (407,127.7) .. (383.5,131.39) ;
\draw    (223,36) .. controls (254,118.25) and (283,132.25) .. (323,36) ;
\draw    (320.26,132.66) .. controls (288.63,50.65) and (259.52,36.87) .. (220.26,133.42) ;
\draw    (325.94,85.55) -- (368.55,85.55) ;
\draw [shift={(370.55,85.55)}, rotate = 180] [color={rgb, 255:red, 0; green, 0; blue, 0 }  ][line width=0.75]    (10.93,-3.29) .. controls (6.95,-1.4) and (3.31,-0.3) .. (0,0) .. controls (3.31,0.3) and (6.95,1.4) .. (10.93,3.29)   ;

\draw (264,160.4) node [anchor=north west][inner sep=0.75pt]    {$I_{2}$};
\draw (387,160.63) node [anchor=north west][inner sep=0.75pt]  [color={rgb, 255:red, 208; green, 2; blue, 27 }  ,opacity=1 ]  {$[ 2,2]$};
\draw (336.11,52.95) node [anchor=north west][inner sep=0.75pt]    {$\pi $};

\end{tikzpicture} \hspace{-5mm}
 & \hspace{-4mm} \begin{tikzpicture}[x=0.75pt,y=0.75pt,yscale=-1,xscale=1,scale=0.4]

\draw    (397,1.7) -- (398.5,161.81) ;
\draw [color={rgb, 255:red, 208; green, 2; blue, 27 }  ,draw opacity=1 ]   (429.55,22.7) .. controls (413.55,25.98) and (394.37,18.69) .. (397.37,48.25) .. controls (400.37,77.81) and (373.44,52.27) .. (358.47,83.23) .. controls (343.5,114.19) and (357.55,111.03) .. (428.55,121.7) ;
\draw    (224.55,45.65) .. controls (240.29,69) and (254.53,90.65) .. (289,89.65) ;
\draw    (289,89.65) .. controls (265.62,93.25) and (244.15,96.2) .. (224.79,134.28) ;
\draw    (297.94,89.7) -- (340.55,89.7) ;
\draw [shift={(342.55,89.7)}, rotate = 180] [color={rgb, 255:red, 0; green, 0; blue, 0 }  ][line width=0.75]    (10.93,-3.29) .. controls (6.95,-1.4) and (3.31,-0.3) .. (0,0) .. controls (3.31,0.3) and (6.95,1.4) .. (10.93,3.29)   ;

\draw (228,158.45) node [anchor=north west][inner sep=0.75pt]    {$II$};
\draw (368,163.24) node [anchor=north west][inner sep=0.75pt]  [color={rgb, 255:red, 208; green, 2; blue, 27 }  ,opacity=1 ]  {$[ 3,1]$};
\draw (308.11,57.1) node [anchor=north west][inner sep=0.75pt]    {$\pi $};

\end{tikzpicture} \hspace{-5mm}
 & \hspace{-3mm}\begin{tikzpicture}[x=0.75pt,y=0.75pt,yscale=-1,xscale=1,scale=0.4]

\draw    (421,6.4) -- (422.5,150.48) ;
\draw [color={rgb, 255:red, 208; green, 2; blue, 27 }  ,draw opacity=1 ]   (370.45,27.72) .. controls (413,39.28) and (418.75,51.84) .. (421.75,78.44) .. controls (424.75,105.04) and (413,123.51) .. (375.5,127.58) ;
\draw    (246.83,20.75) .. controls (272.61,90.3) and (296.73,102.14) .. (330,20.75) ;
\draw    (327.72,133.77) .. controls (301.41,64.42) and (277.2,52.77) .. (244.55,134.42) ;
\draw    (321.94,77.4) -- (364.55,77.4) ;
\draw [shift={(366.55,77.4)}, rotate = 180] [color={rgb, 255:red, 0; green, 0; blue, 0 }  ][line width=0.75]    (10.93,-3.29) .. controls (6.95,-1.4) and (3.31,-0.3) .. (0,0) .. controls (3.31,0.3) and (6.95,1.4) .. (10.93,3.29)   ;

\draw (257,150.4) node [anchor=north west][inner sep=0.75pt]    {$III$};
\draw (377,148.01) node [anchor=north west][inner sep=0.75pt]  [color={rgb, 255:red, 208; green, 2; blue, 27 }  ,opacity=1 ]  {$[ 4]$};
\draw (332.11,44.8) node [anchor=north west][inner sep=0.75pt]    {$\pi $};

\end{tikzpicture} \hspace{-3mm}
 \\

 \hline

\end{tabular}
}
\caption{Singular fibers of the elliptic fibration induced by smooth fibers of the conic bundle.}\label{table:smooth fibers}
\end{table}



\begin{table}[H] 

\setlength{\tabcolsep}{5mm} 
\def\arraystretch{1.25} 
\centering
\makebox[\linewidth]{
\begin{tabular}{|c|c|c|c|c|}
\hline
\hspace{-3mm}\begin{tikzpicture}[x=0.75pt,y=0.75pt,yscale=-1,xscale=1,scale=0.4]

\draw [color={rgb, 255:red, 65; green, 117; blue, 5 }  ,draw opacity=1 ]   (429.52,7.7) -- (402.33,84.69) ;
\draw [color={rgb, 255:red, 208; green, 2; blue, 27 }  ,draw opacity=1 ]   (373.45,19.14) .. controls (406.08,14.81) and (425.57,12.53) .. (428.83,38.51) .. controls (432.09,64.48) and (402.81,42.04) .. (386.53,69.24) .. controls (370.24,96.45) and (428.43,85.09) .. (429.52,103.85) .. controls (430.6,122.61) and (404.5,113.23) .. (378.94,116.65) ;
\draw [color={rgb, 255:red, 65; green, 117; blue, 5 }  ,draw opacity=1 ]   (237.88,33.05) .. controls (264.27,105.4) and (288.95,117.72) .. (323,33.05) ;
\draw [color={rgb, 255:red, 74; green, 144; blue, 226 }  ,draw opacity=1 ]   (320.66,118.07) .. controls (293.74,45.93) and (268.97,33.81) .. (235.55,118.75) ;
\draw [color={rgb, 255:red, 74; green, 144; blue, 226 }  ,draw opacity=1 ]   (402.33,61.82) -- (425.17,129.75) ;
\draw    (322.94,77.4) -- (365.55,77.4) ;
\draw [shift={(367.55,77.4)}, rotate = 180] [color={rgb, 255:red, 0; green, 0; blue, 0 }  ][line width=0.75]    (10.93,-3.29) .. controls (6.95,-1.4) and (3.31,-0.3) .. (0,0) .. controls (3.31,0.3) and (6.95,1.4) .. (10.93,3.29)   ;

\draw (274,132.45) node [anchor=north west][inner sep=0.75pt]    {$I_{2} \quad [1,1],[1,1]$ };
\draw (333.11,44.8) node [anchor=north west][inner sep=0.75pt]    {$\pi $};

\end{tikzpicture}
\hspace{-4mm} & \hspace{-4mm} \begin{tikzpicture}[x=0.75pt,y=0.75pt,yscale=-1,xscale=1,scale=0.4]

\draw [color={rgb, 255:red, 65; green, 117; blue, 5 }  ,draw opacity=1 ]   (445,6.6) -- (400,81.3) ;
\draw [color={rgb, 255:red, 208; green, 2; blue, 27 }  ,draw opacity=1 ]   (373.45,10.24) .. controls (403.45,5.54) and (421.37,3.07) .. (424.37,31.23) .. controls (427.37,59.39) and (400.44,35.06) .. (385.47,64.56) .. controls (370.5,94.05) and (424,81.73) .. (425,102.07) .. controls (426,122.41) and (402,112.24) .. (378.5,115.95) ;
\draw [color={rgb, 255:red, 74; green, 144; blue, 226 }  ,draw opacity=1 ]   (320.26,92.66) .. controls (288.63,10.65) and (259.52,-3.13) .. (220.26,93.42) ;
\draw [color={rgb, 255:red, 74; green, 144; blue, 226 }  ,draw opacity=1 ]   (400,56.5) -- (421,130.15) ;
\draw [color={rgb, 255:red, 65; green, 117; blue, 5 }  ,draw opacity=1 ]   (306,32.3) -- (248,106.7) ;
\draw [color={rgb, 255:red, 65; green, 117; blue, 5 }  ,draw opacity=1 ]   (230,38.1) -- (295,107.3) ;
\draw    (323.94,69.4) -- (366.55,69.4) ;
\draw [shift={(368.55,69.4)}, rotate = 180] [color={rgb, 255:red, 0; green, 0; blue, 0 }  ][line width=0.75]    (10.93,-3.29) .. controls (6.95,-1.4) and (3.31,-0.3) .. (0,0) .. controls (3.31,0.3) and (6.95,1.4) .. (10.93,3.29)   ;

\draw (264,125.4) node [anchor=north west][inner sep=0.75pt]    {$I_{3}\quad [2],[1,1]$};
\draw (334.11,36.8) node [anchor=north west][inner sep=0.75pt]    {$\pi $};

\end{tikzpicture} \hspace{-5mm}
 & \hspace{-4mm} \begin{tikzpicture}[x=0.75pt,y=0.75pt,yscale=-1,xscale=1,scale=0.4]

\draw [color={rgb, 255:red, 65; green, 117; blue, 5 }  ,draw opacity=1 ]   (442,3.9) -- (397,91.95) ;
\draw [color={rgb, 255:red, 208; green, 2; blue, 27 }  ,draw opacity=1 ]   (370.45,8.19) .. controls (400.45,2.66) and (418.37,-0.26) .. (421.37,32.94) .. controls (424.37,66.13) and (397.44,37.45) .. (382.47,72.22) .. controls (367.5,106.98) and (421,92.46) .. (422,116.43) .. controls (423,140.41) and (399,128.42) .. (375.5,132.8) ;
\draw [color={rgb, 255:red, 74; green, 144; blue, 226 }  ,draw opacity=1 ]   (397,62.73) -- (436,142.25) ;
\draw [color={rgb, 255:red, 65; green, 117; blue, 5 }  ,draw opacity=1 ]   (302,54.3) -- (244,128.7) ;
\draw [color={rgb, 255:red, 65; green, 117; blue, 5 }  ,draw opacity=1 ]   (226,60.1) -- (291,129.3) ;
\draw [color={rgb, 255:red, 74; green, 144; blue, 226 }  ,draw opacity=1 ]   (253,20.25) -- (302,85.25) ;
\draw [color={rgb, 255:red, 74; green, 144; blue, 226 }  ,draw opacity=1 ]   (274,18.25) -- (220,84.25) ;
\draw    (315.94,83.75) -- (358.55,83.75) ;
\draw [shift={(360.55,83.75)}, rotate = 180] [color={rgb, 255:red, 0; green, 0; blue, 0 }  ][line width=0.75]    (10.93,-3.29) .. controls (6.95,-1.4) and (3.31,-0.3) .. (0,0) .. controls (3.31,0.3) and (6.95,1.4) .. (10.93,3.29)   ;

\draw (258,142.4) node [anchor=north west][inner sep=0.75pt]    {$I_{4}\quad [2],[2]$};
\draw (326.11,51.15) node [anchor=north west][inner sep=0.75pt]    {$\pi $};

\end{tikzpicture} \hspace{-5mm}
 & \hspace{-4mm} \begin{tikzpicture}[x=0.75pt,y=0.75pt,yscale=-1,xscale=1,scale=0.4]

\draw [color={rgb, 255:red, 65; green, 117; blue, 5 }  ,draw opacity=1 ]   (228.94,18) .. controls (255.93,79.46) and (281.17,89.92) .. (316,18) ;
\draw [color={rgb, 255:red, 74; green, 144; blue, 226 }  ,draw opacity=1 ]   (313.61,117.88) .. controls (286.07,56.59) and (260.73,46.3) .. (226.55,118.45) ;
\draw [color={rgb, 255:red, 65; green, 117; blue, 5 }  ,draw opacity=1 ]   (413,5.05) -- (379.45,78.92) ;
\draw [color={rgb, 255:red, 208; green, 2; blue, 27 }  ,draw opacity=1 ]   (361.45,17.37) .. controls (391.45,12.7) and (413.45,21.41) .. (412.37,38.23) .. controls (411.28,55.05) and (388.44,42.03) .. (373.47,71.33) .. controls (358.5,100.64) and (405.45,117.77) .. (440.45,124.77) ;
\draw [color={rgb, 255:red, 74; green, 144; blue, 226 }  ,draw opacity=1 ]   (386.45,31.51) -- (409,136.5) ;
\draw    (305.94,72.1) -- (348.55,72.1) ;
\draw [shift={(350.55,72.1)}, rotate = 180] [color={rgb, 255:red, 0; green, 0; blue, 0 }  ][line width=0.75]    (10.93,-3.29) .. controls (6.95,-1.4) and (3.31,-0.3) .. (0,0) .. controls (3.31,0.3) and (6.95,1.4) .. (10.93,3.29)   ;

\draw (241,150) node [anchor=north west][inner sep=0.75pt]    {$III \quad [1,\textcolor{verde}{\mathbf{1}}],[\textcolor{azul}{\mathbf{1}},1]$};
\draw (316.11,39.5) node [anchor=north west][inner sep=0.75pt]    {$\pi $};

\end{tikzpicture} \hspace{-5mm}
 &\hspace{-4mm} \begin{tikzpicture}[x=0.75pt,y=0.75pt,yscale=-1,xscale=1,scale=0.4]

\draw [color={rgb, 255:red, 65; green, 117; blue, 5 }  ,draw opacity=1 ]   (202.45,28.24) .. controls (230.4,94.69) and (280.49,103.21) .. (316.55,25.45) ;
\draw [color={rgb, 255:red, 65; green, 117; blue, 5 }  ,draw opacity=1 ]   (415.1,2.45) -- (369.45,93.73) ;
\draw [color={rgb, 255:red, 208; green, 2; blue, 27 }  ,draw opacity=1 ]   (451.24,24.98) .. controls (391.32,18.64) and (374.21,42.58) .. (380.86,72.15) .. controls (387.52,101.72) and (445.53,94.68) .. (464.55,56.66) ;
\draw [color={rgb, 255:red, 74; green, 144; blue, 226 }  ,draw opacity=1 ]   (376.11,50.78) -- (397.55,145.9) ;
\draw [color={rgb, 255:red, 74; green, 144; blue, 226 }  ,draw opacity=1 ]   (224.99,39.55) -- (279.09,111.45) ;
\draw [color={rgb, 255:red, 74; green, 144; blue, 226 }  ,draw opacity=1 ]   (284.5,37.93) -- (234.91,111.45) ;
\draw    (316.94,77.2) -- (359.55,77.2) ;
\draw [shift={(361.55,77.2)}, rotate = 180] [color={rgb, 255:red, 0; green, 0; blue, 0 }  ][line width=0.75]    (10.93,-3.29) .. controls (6.95,-1.4) and (3.31,-0.3) .. (0,0) .. controls (3.31,0.3) and (6.95,1.4) .. (10.93,3.29)   ;

\draw (242,150) node [anchor=north west][inner sep=0.75pt]    {$IV \quad [\textcolor{azul}{\mathbf{2}}],[\textcolor{verde}
{\mathbf{1}},1]$};

\draw (327.11,44.6) node [anchor=north west][inner sep=0.75pt]    {$\pi $};

\end{tikzpicture} \hspace{-5mm}
 \\ \hline
\end{tabular}
}
\caption{Singular fibers of the elliptic fibration induced by singular fibers of the conic bundle.}\label{table:sing_fibers}
\end{table}

\end{proposition}

\begin{proof}
Let $F$ be a smooth fiber of $f$. By the previous construction, the branching locus of $\pi:X\rightarrow Z$ is a smooth irreducible curve $C$ (in red in Table \ref{table:smooth fibers}) and by Remark \ref{rem_singular_fibers}, $F$ meets $C$ in 4 points (with multiplicity). We study each case separately.

If $C$ meets $F$ in 4 distinct point, then the corresponding fiber of $\mathcal E$ is a double cover of $\mathbf P^1$ with 4 ramification points.
In Table \ref{table:smooth fibers} we show the possible multiplicities of the points.
Then by Riemann-Hurwitz formula the fiber of $\mathcal E$ is a smooth curve of genus 1, i.e. a curve of type $I_0$.
If $C$ meets $F$ in 2 simple points and a double point $p$, the preimage of $p$ is a nodal point of the fiber, thus the fiber is of type $I_1$.
Similarly one obtains a fiber of type $I_2$ when $C\cap F$ consists of two double points.
If $C$ meets $F$ in a triple point and a simple one, then the fiber of $\mathcal E$ has a singular point which is a cusp: the triple point in the double cover gives a singular point of the fiber where the equation is locally given by $y^2=x^3$, thus a cusp.
Similarly, if $C$ meets $F$ in a single point with multiplicity 4 the induced fiber of $\mathcal E$ is of type $III$.

Now let $F$ be a singular fiber of $f$. According to Remark \ref{rem_singular_fibers},  singular fibers of a conic bundle on $Z$ are union of two $(-1)$-curves $F_1$ and $F_2$ that intersect at one point $P=F_1\cap F_2$. Furthermore, the branch curve $C$ of the $2$-cover $X\to Z$ (in red in Table \ref{table:sing_fibers}) intersects each rational curve $F_i$ at two points (with multiplicity). In Table \ref{table:sing_fibers} we show the possible multiplicities of the points.
If $C$ meets each $F_i$ in two simple points, distinct from $P$, this defines a fiber of type $I_2$ of $\mathcal E$, since the preimage of each $F_1$ defines a rational curve and $\pi^{-1}(P)$ consists of two points.
If $C$ meets $F_2$ in two simple points and $F_1$ in a double point, distinct from $P$, this defines a singular fiber of $\mathcal E$ of type $I_3$: the double cover of $F_1$ gives two components of the fiber $I_3$, while $F_2$ contributes with one component and $\pi^{-1}(P)$ consist of two points.
Similarly, if $C$ meets each $F_i$ in a double point, distinct from $P$, the fiber is of type $I_4$. The last case to study is when $C$ meets $F_1$ in two simple points $F_2$ in $P$ with multiplicity 2. In this case, the fiber $\pi^{-1}(F)$ has three components and they all meet in $\pi^{-1}(P)$. Three concurrent rational curves form  a fiber of type $IV$. 
\end{proof}

Now we want to classify which types of fibers are compatible in each case with the classification given in Proposition \ref{Classes of conic bundles on Z}. Note that in the case $Z\simeq \mathbf{P}^1\times \mathbf{P}^1$ all types of singular fibers in Proposition \ref{theoretical fibers} (1) can be realized, as Example \ref{example:P1xP1} shows. Now, given a del Pezzo surface $Z_d$ (see Convention \ref{convention:Z_d}), Theorem \ref{thm_main} establishes a correspondence between conic bundles on $Z_d$ and elliptic fibrations on the K3 surface $X$.  Proposition \ref{Classes of conic bundles on Z} then classifies conic bundles on $Z_d$.
Given the potential fibers outlined in Proposition \ref{theoretical fibers}, our goal is to determine which ones are admissible for each $Z_d$ and to establish their connection with the geometry of the sextic curve $\Gamma_d$ introduced in Proposition \ref{prop: K3 of strictly elliptic type}.

\begin{proposition}\label{prop: admissible fibers}
   Let $(X,\iota)$ be a pair of strictly elliptic type and $Z_d=X\slash \langle \iota \rangle$ be the quotient smooth del Pezzo surface as in Proposition \ref{prop: K3 of strictly elliptic type}. If $d\leq 5$, all types of singular fibers in Proposition \ref{theoretical fibers} are admissible. 
    If $d=6,7$ (resp. 8) then the fiber $I_4$ and $IV$ (resp. $I_3$, $I_4$ and $IV$) are not admissible. 
    
    In other words, the admissible fibers for the elliptic fibration $\mathcal{E}:X\to \mathbf{P}^1$ are described in the following table

    \begin{table}[H]
        \centering
        \begin{tabular}{c|l}
            $d$& \textup{Singular fibers}\\
            \hline
            $8$ &$I_0, I_1, I_2, II, III$\\
            $7,6$&  $I_0, I_1, I_2, I_3, II, III$\\
            $\leq5$ & $I_0, I_1, I_2, I_3, I_4, II, III, IV$\\
        \end{tabular}
        \caption{Admissible singular fibers of the elliptic fibration $\mathcal{E}:X\to \mathbf{P}^1$ induced by a conic bundle $f:Z_d\to \mathbf{P}^1$.}
        \label{tab_fibras_admisibles}
    \end{table}
\end{proposition}
\begin{proof}
    Let $d=8$. By Proposition \ref{Classes of conic bundles on Z}, the only conic class on $Z_8$ is $D=L-E_1$, where $E_1$ is the exceptional divisor of the blow up. 
    In this case the conic bundle has no singular fibers, thus $\Gamma_d$ meets the fiber in the nodal point and four more points.
    According to possibilities given in Table \ref{table:smooth fibers}, one can have singular fiber of the elliptic fibrations of types $I_0$, $I_1$, $I_2$, $II$ or $III$. 

If $d=7$, the curve $\Gamma_d$ is a sextic curve with 2 nodes and by Proposition \ref{Classes of conic bundles on Z}, canonical classes on $Z_7$ are $D=L-E_i$, where $E_i$ is the exceptional divisor over one of the two nodal points $p_1,p_2$. We take $i=1$ without loss of generality.
If the conic bundle has no singular fibers, then as before the possible fibers of the elliptic fibration are of types $I_0$, $I_1$, $I_2$, $II$ or $III$.
If $F=F_1+F_2$ reducible, which corresponds to the case when $L$ passes through the other nodal point $p_2$, we observe that each component $F_1$ and $F_2$ meets $\Gamma_d$ in 2 points. If they are distinct for both $F_1,F_2$, one has again a fiber of type $I_2$, while is one of the two components meets $\Gamma_d$ with multiplicity 2, one obtains a fiber of type $I_3$ (see case $[1,1],[1,1]$ and case $[2],[1,1]$ of Table \ref{table:sing_fibers}).
Observe that the strict transform always meet one of the two components in two distinct points (coming from the blow up of the nodal points). Thus the only cases of Table \ref{table:sing_fibers} that appear are $[1,1],[1,1]$ and $[2],[1,1]$.
The case $d=6$ is analogous.
 
When $d=5$ according to Proposition \ref{Classes of conic bundles on Z}, one can have the conic class $D=2L-\sum_{i=1}^4 E_i$. In this case we distinguish if conics of the bundle are irreducible or reducible. 
The first case will give fibers of type $I_0, I_1, I_2, II$ and $III$ as in Table \ref{table:smooth fibers}. 
On the other hand, in the case of conics reduced as the union of two lines $L_1\cup L_2$, one observes that this allows the cases of Table \ref{table:sing_fibers}: each line passes through 2 of the 5 nodal points of the sextic $\Gamma_5$, thus meeting $\Gamma_5$ in two more points (with multiplicity). According to the distribution of these points, one obtain all cases of Table \ref{table:sing_fibers}.

The remaining cases with $d\leq 4$ can be treated in the same way, and thus admit all types of fibers listed in Proposition \ref{theoretical fibers}.
\end{proof}

\begin{remark}
    The proof of the previous result, along with Table \ref{tab_canonical_bundles}, not only allows for the determination of admissible singular fibers in the induced elliptic fibration but also the complete configuration of these fibers.
\end{remark}

\bibliographystyle{alpha} 
\bibliography{bibliography}

\end{document}